\documentclass[10pt,a4paper]{article}
\usepackage{amsmath}
\usepackage{amsfonts}
\usepackage{amsmath, amssymb}
\usepackage{makeidx}
\makeindex

\newtheorem{theorem}{Theorem}[section]
\newtheorem{lemma}{Lemma}[section]

\newtheorem{proposition}{Proposition}[section]

\newtheorem{definition}{Definition}[section]

\newtheorem{remark}[theorem]{Remark}

\newcommand{\R}{{{\mathbb R}}}
\newcommand{\N}{{{\mathbb N}}}

\numberwithin{equation}{section}

\begin{document}
\author{Yanjin Wang \thanks{Corresponding email: wangyj@ms.u-tokyo.ac.jp}
\\{\small \textit{Graduate School
of Mathematical Sciences, University of Tokyo}}
\\{\small\textit{3-8-1 Komaba, Meguro, Tokyo, 153-8914, Japan}}}

\title{Nonexistence of global solutions of a class of coupled nonlinear Klein-Gordon equations with nonnegative potentials and arbitrary initial energy}
\date{}
\maketitle

\begin{abstract}
In the paper we consider the nonexistence of global solutions of the
Cauchy problem for coupled Klein-Gordon equations of the form
\begin{eqnarray*}
\left\{\begin{array}{l}
u_{tt}-\Delta u+m_1^2 u+K_1(x)u=a_1|v|^{q+1}|u|^{p-1}u\\
v_{tt}-\Delta v+m_2^2 u+K_2(x)v=a_2|u|^{p+1}|v|^{q-1}v\\
u(0,x)=u_0; u_t(0,x)=u_1(x)\\
v(0,x)=v_0; v_t(0,x)=v_1(x)
\end{array}
\right.
\end{eqnarray*}
on $\R\times\R^n$.

Firstly for some special cases of $n=2,3$, we prove the existence of
ground state of the corresponding Lagrange-Euler equations of the
above equations. Then we establish a blow up result with low initial
energy, which leads to instability of standing waves of the system
above. Moreover as a byproduct we also discuss the global existence.
Next based on concavity method we prove the blow up result for
the system with non-positive initial energy in the general case:
$n\geq 1$. Finally when the initial energy is given arbitrarily positive,
we show that if the initial datum satisfies some conditions,
the corresponding solution blows up in a finite time.


\end{abstract}

\section{Introduction}
The motion of charged measons in an electromagnetic field can be
described by the following coupled Klein-Gordon equations:
\begin{eqnarray}
\left\{
\begin{array}{l}
u_{tt}-\Delta u+\alpha^2 u+g^2v^2u=0\\
v_{tt}-\Delta v+\beta^2 v+h^2u^2v=0
\end{array}
\right.\label{KGE}
\end{eqnarray}
where $\Delta$ is Laplacian operator on $\R^n$, $\alpha$ and $\beta$
are non zero real constants. The system was firstly introduced by I.
Segal \cite{00Segal}. A lot of authors have discussed this mixed
system; see for example \cite{01Jorgen}, \cite{06LT},
\cite{02Makhandov}, \cite{03Medeiros}. And a sharp condition for
global existence and blowing up has been given by
Zhang \cite{031Zhang} for the mixed problem (\ref{KGE}), where the
blow up result is given under the condition that the initial energy is
below the energy wall. The system was generalized by Miranda and
Mederiros \cite{04Medeiros}, \cite{05Miranda} for the case where the
nonlinear terms are of the form $|v|^{\rho+1}|u|^{\rho-1}u$ and
$|u|^{\rho+1}|v|^{\rho-1}v$ with some $\rho$. Li and
Tsai \cite{07LT2} recently considered a class of nonlinear term which
includes the above generalized nonlinear terms on a bounded domain
of $\R^n$. And some other nonlinear term was considered by Delort,
Fang and Xue \cite{08DelartFX} on $\R^2$.

In this paper we are interested in the initial boundary value
problem for the coupled Klein-Gordon equations with nonnegative
potentials of the form
\begin{eqnarray}
\left\{\begin{array}{l}
u_{tt}-\Delta u+m_1^2 u+K_1(x)u=a_1|v|^{q+1}|u|^{p-1}u\\
v_{tt}-\Delta v+m_2^2 v+K_2(x)v=a_2|u|^{p+1}|v|^{q-1}v\\
u(0,x)=u_0; u_t(0,x)=u_1(x), x\in\R^n\\
v(0,x)=v_0; v_t(0,x)=v_1(x), x\in\R^n
\end{array}
\right.\label{CKGE_P}
\end{eqnarray}
where the parameters $a_1$ and $a_2$ are positive constants, the
masses are nonzero, $m_1\ne 0$ and $m_2\ne 0$.

 In the paper we make the following restriction on the
real numbers $p>1$ and $q>1$:

If $n=1,2$,
\begin{eqnarray}
1<p, q<\infty;\label{pq-cond1}
\end{eqnarray}
And if $n\geq 3$, then
\begin{eqnarray}
q<p+1<\frac{n+2}{n-2} \mathrm{\ or\ }
p+1<q<\frac{n+2}{n-2},\label{pq-cond2}
\end{eqnarray}
and
\begin{eqnarray}
p<q+1<\frac{n+2}{n-2} \mathrm{\ or\ }
q+1<p<\frac{n+2}{n-2}.\label{pq-cond3}
\end{eqnarray}

Throughout the paper we will assume that $K_i(x)$ satisfies
\begin{eqnarray}
K_i(x)\geq 0 \mathrm{\ \ (}\forall x\in\R^n\mathrm{)}
\end{eqnarray}
for $i=1,2$.

Before describing our results, we first recall the existing
results about the Cauchy problem for the single Klein-Gordon
equation
\begin{eqnarray}
u_{tt}-\Delta u+m^2 u+K(x) u=f(u) \mathrm{\ on\
}[0,\infty)\times\R^n,\label{SKGE}
\end{eqnarray}
where $K(x)\geq 0$.

When $K(x)\equiv 0$, there are numerous results about the existence
and blow up of solutions for the equation (\ref{SKGE}). It is well
known that the solution blows up in a finite time when the initial
energy is negative. Here we refer to
\cite{11Ginibre}, \cite{12Ginibre}, \cite{10Levine}, \cite{09Strauss}. The
instability of standing wave of the equation (\ref{SKGE}) was
studied in
\cite{14Berestycki}, \cite{141Berestycki}, \cite{26Liu}, \cite{13Ohta}.
Based on the results in
\cite{14Berestycki}, \cite{141Berestycki}, Zhang \cite{15Zhang}
established a sharp condition of global existence and blow up for
the eqution (\ref{SKGE}). Recently the author \cite{16wang} has
proposed a sufficient condition of the initial datum with
arbitrarily positive initial energy such that the corresponding
solution of the equation (\ref{SKGE}) blows up in a finite time.

As for the case $K(x)\not\equiv 0$, when the nonlinear term
vanishes it has been shown that the equation (\ref{SKGE}) has time
periodic and spatially localized solutions in \cite{17Rose}. And
 Soffer and Weinsterin \cite{18Soffer} considered a class of
nonlinear Klein-Gordon equations (\ref{SKGE}) with nonlinear term
$f(u)$, which is real-valued, smooth in a neighborhood of $u=0$ and
has an expansion $f(u)=u^3+O(u^4)$ on $\R^3$. Recently Gan and
Zhang \cite{19Gan} considered standing waves for the equation
(\ref{SKGE}) with $f(u)=|u|^{p-1}u$.

Now we return to the coupled Klein-Gordon equations (\ref{CKGE_P}).
As we know, until now there is no result for the system
(\ref{CKGE_P}). In the paper we are concerned with the nonexistence of
global solutions of the system (\ref{CKGE_P}). We first establish
the existence of a unique local weak solution of the equations
(\ref{CKGE_P}) by applying the Banach constraction mapping
principle.

Based on the local existence theorem, our first purpose of the
present paper is to establish a blowing up result by using the
ground state solution. This blowing up result leads to instability
of standing wave for the equation (\ref{CKGE_P}). The proof is done
by first showing the existence of the ground state solution of the
corresponding Lagrange-Euler equations by variational method which
was firstly introduced in \cite{14Berestycki}, \cite{141Berestycki},
and then by discussing a blow up result with low initial energy
based on a potential well argument and concavity method, which is
originated by Payne and Sattinger \cite{24Payne} and Levine
\cite{10Levine}, \cite{25Levine}, respectively. Because of the
restriction of the embedding theorem of $H^1_K\hookrightarrow L^r$
($\displaystyle\displaystyle 2<r<\frac{2n}{n-2}$) the blow up result
will be established only on $\R^n$ ($n=2,3$). As a byproduct we,
however, also establish the global existence of solutions of the
system (\ref{CKGE_P}) when $n=2,3$.

Our next purpose is to show the blow up result when the initial
energy is non-positive by a concavity argument.

The final purpose is to construct sufficient conditions of the
initial datum such that the corresponding solution blows up in
a finite time with arbitrarily positive initial energy, that is, we
show that there exists a finite time $T_{\max}$ with $E(0)>0$ such
that
\begin{eqnarray*}
\lim_{t\rightarrow T_{\max}^-}(\|u(t)\|^2_2+\|v(t)\|^2_2)=\infty.
\end{eqnarray*}
To the best of our knowledge, this is the first blowing up result
for the coupled Klein-Gordon equations with arbitrarily positive
initial energy.

The paper is organized as follows. In Section 2 we state the local
existence of solutions by a fixed point argument. In Section 3 we
study the existence of the standing wave of the system
(\ref{CKGE_P}) with the ground state by using variational method. In
section 4, based on the result obtained in Section 3 we establish
blowing up result for the system (\ref{CKGE_P}) on $\R^n$ (n=2,3),
which will lead to the instability of standing waves. In Section 5,
using concavity argument we show a blowing up result when the
initial energy is negative. In the last section, we establish some
sufficient conditions of initial datum with arbitrarily initial
energy such that the corresponding solution blows up in a finite time.

\section{Local existence}

In the paper we will work in the energy space:
\begin{eqnarray*}
H_{0, K_i}^1=\left\{u\in H_0^1(\R^n):\int K_i(x)|u(x)|^2dx<\infty
\right\}
\end{eqnarray*}
 with the following norm:
\begin{eqnarray*}
\|u\|_{H^1_{0, K_i}}=\|\nabla u\|^2+m_i^2\|u\|^2+\int_{\R^n}
K_i(x)|u(x)|^2dx
\end{eqnarray*}
where $K_i(x)\geq 0$ ($\forall x\in\R^n$) for $i=1,2$.

And for simplicity we denote $\displaystyle\int_{\R^n}dx$ by
$\displaystyle\int dx$. The notation $t\rightarrow T^-$ means that
$t\rightarrow T$ and $t<T$.

Firstly we rewrite the coupled Klein-Gordon equations (\ref{CKGE_P})
in the following equivalent form
\begin{eqnarray}
\left\{\begin{array}{l}
\alpha u_{tt}-\alpha \Delta u+\alpha m_1^2 u+\alpha K_1(x)u=a_2^\prime(p+1)|v|^{q+1}|u|^{p-1}u\\
v_{tt}-\Delta v+m_2^2 u+K_2(x)v=a_2^\prime(q+1)|u|^{p+1}|v|^{q-1}v\\
u(0,x)=u_0; u_t(0,x)=u_1(x)\\
v(0,x)=v_0; v_t(0,x)=v_1(x)
\end{array}
\right.\label{CKGE_P1}
\end{eqnarray}
on $\R\times\R^n$, where
$\displaystyle\alpha=\frac{a_2(p+1)}{a_1(q+1)}$ and $\displaystyle
a_2^\prime=\frac{a_2}{q+1}$.

\begin{definition} A function $(u, v)$ is said to be a solution of
the system (\ref{CKGE_P1}), if it satisfies that $u \in C^0([0,
T), H_{0, K_1}^1 (\R^n)) \cap C^1([0, T), L^2(\R^n))$ and $v\in
C^0([0, T), H_{0, K_2}^1(\R^n))\cap C^1([0, T), L^2(\R^n))$,
\begin{eqnarray*}
\alpha\left[\int u_{tt}w_1dx+\int \nabla u(t) \nabla w_1dx+m_1^2\int
u(t)w_1dx+\int
K_1(x)u(t)w_1dx\right]\\
=a_2^\prime(p+1)\int|v(t)|^{q+1}|u(t)|^{p-1}u(t)w_1dx,\\
\int v_{tt}w_2dx+\int \nabla v(t) \nabla w_2dx+m_2^2\int
v(t)w_2dx+\int
K_2(x)v(t)w_2dx\\
=a_2^\prime(q+1)\int|u(t)|^{p+1}|v(t)|^{q-1}v(t)w_2dx,
\end{eqnarray*}
and
\begin{eqnarray*}
\left\{
\begin{array}{l}
u(0,x)=u_0(x), u_t(0,x)=u_1(x)\\
v(0,x)=v_0(x), v_t(0,x)=v_1(x)
\end{array}
 \right.
\end{eqnarray*}
for $(w_1,w_2)\in H^1_{0, K_1}\times H^1_{0, K_2}$, $x\in\R^n$ and
$t\in[0,T)$.
\end{definition}

Now we state the local existence theorem of the initial boundary
value problem for the equivalent system (\ref{CKGE_P1}).
\begin{theorem} Assume that $p$ and $q$ satisfy the conditions
(\ref{pq-cond1}), (\ref{pq-cond2}) and (\ref{pq-cond3}). Let $(u_0,
v_0)\in H_{0, K_1}^1\times H_{0, K_2}^1$ and $(u_1, v_1)\in
L^2\times L^2$. Then there exists a unique solution $(u(t,x),
v(t,x))$ of the equations (\ref{CKGE_P1}) on a maximal time interval
$[0, T_{\max})$ for some $T_{\max}\in(0,\infty]$ such that $(u,
v)\in C^0([0, T_{\max}); H^1_{0, K_1}(\R^n))\times C^0([0,
T_{\max}); H^1_{0, K_2}(\R^n))$.

Furthermore, we have the following alternatives:
\begin{eqnarray}
T_{\max}=\infty;
\end{eqnarray}
or
\begin{eqnarray}
T_{\max}<\infty \mathrm{\ and\ } \lim_{t\rightarrow
T_{\max}^-}(\alpha\|u(t)\|_2^2+\|v(t)\|_2^2)=\infty.
\end{eqnarray}

Moreover, the local solution $(u,v)$ satisfies the following
conservation law of energy:
\begin{eqnarray}
E(t)=E(0)\label{E_conservation}
\end{eqnarray}
for every $t\in[0, T_{\max})$, where
\begin{eqnarray}
E(t)&=&\frac{1}{2}\int(\alpha|u_t|^2+|v_t|^2)dx+\frac{1}{2}\left[\int
\alpha\left(|\nabla u(t)|^2+m_1^2|u(t)^2|+K_1(x)|u(t)|^2\right)dx\right.\nonumber\\
&&+\left.\int\left(|\nabla
v(t)|^2+m_2^2|v(t)^2|+K_2(x)|v(t)|^2\right)dx\right]
-a_2^\prime\int|u(t)|^{p+1}|v(t)|^{q+1}dx. \label{E_def}
\end{eqnarray}
\end{theorem}

In order to prove the above theorem, we consider the following
scalar equation
\begin{eqnarray}
\left\{
\begin{array}{l}
w_{tt}-\Delta w+m^2w+K(x)w=f(t,x)\\
w(0,x)=w_0, w_t(0,x)=w_1(x)
\end{array}
\right.\label{CKGE_P2}
\end{eqnarray}
where $K(x)\geq 0$ and $m\ne 0$.

\begin{theorem}
Assume that $f(t,x)$ is a Lipschitz function with respect to $x$. If
$(w_0,w_1)\in H^1_{0, K}(\R^n)\times L^2(\R^n)$, then there exists
a unique function $w\in C^1([0, T], H^1_{0, K})$ satisfying the
equation (\ref{CKGE_P2}) for $T>0$.
\end{theorem}

The proof of this theorem follows the argument by
Haraux \cite{20Haraux}, Lions and Magens \cite{21Lions}. We here omit
it.

We next give two estimates on the nonlinear term of the system
(\ref{CKGE_P1}).

\begin{lemma}
Assume that $p$ and $q$ satisfy the conditions (\ref{pq-cond1}),
(\ref{pq-cond2}) and (\ref{pq-cond3}). Then we have the following
estimates:
\begin{eqnarray}
\||v|^{q+1}|u|^{p-1}uw\|_{1}\leq\|\nabla v\|^{q+1}\|\nabla u\|^p\|\nabla w\|,\\
\||u|^{p+1}|v|^{q-1}vw\|_1\leq\|\nabla u\|^{p+1}\|\nabla
u\|^q\|\nabla w\|.
\end{eqnarray}
\end{lemma}
\textbf{Proof.} Firstly we consider the case $n\geq 3$ and $q+1>p$.
By H\"older inequality and Sobolev inequality we have
\begin{eqnarray*}
\||v|^{q+1}|u|^{p-1}u
w\|_1\leq\|vu\|^p_{r_1}\|v\|^{q+1-p}_{r_2}\|w\|_{r_3}\\
\leq \|\nabla v\|^{q+1}\|\nabla u\|^p\|\nabla w\|,
\end{eqnarray*}
where $r_1$, $r_2$ and $r_3$ satisfy the following conditions:
\begin{eqnarray*}
&&\frac{p}{r_1}+\frac{q+1-p}{r_2}+\frac{1}{r_3}=1,\\
&&2<r_1<\frac{2n}{n-2},\\
&&2<r_2<\frac{2n}{n-2},\\
&&2<r_3<\frac{2n}{n-2}.
\end{eqnarray*}
From the inequalities above on $r_1$, $r_2$ and $r_3$, we get
that $\displaystyle q<\frac{4}{n-2}$.

For the other cases, using the same argument as above we can obtain the
desired result.
\begin{flushright}
$\Box$
\end{flushright}

\textbf{Proof of Theorem 2.1.} The proof relies on the Banach
contraction mapping principle. For $T>0$ and $R>0$ we define the
following space
\begin{eqnarray*}
X_{T,R}\equiv\{(u(t),v(t));u(t)\in C^0([0, T];
H^1_{0,K_1}(\R^n))\cap L^2([0,T];H^1_{0,K_1}),\nonumber\\
v(t)\in C^0([0, T];
H^1_{0,K_2}(\R^n))\cap L^2([0,T];H^1_{0,K_2}),\nonumber\\
(u_t,v_t)\in C^1([0, T]; L^2(\R^n))\times C^1([0, T]; L^2(\R^n)),\nonumber\\
 e(u,v)\leq R^2 \mathrm{\ for\ all\ }t\in[0, T],
 \nonumber\\
 (u(0),v(0))=(u_0,v_0)\mathrm{\ and\
 }(u_t(0),v_t(0))=(u_1,v_1)\}\nonumber\\
\end{eqnarray*}
where
$$\displaystyle e(u,v)=\max_{t\in[0,
T]}\{\alpha(\|u(t)\|_{H^1_{0,
K_1}}^2+\|u_t(t)\|_2^2)+\|v(t)\|_{H^1_{0,
K_2}}^2+\|v_t(t)\|_2^2\}.$$

We now define a nonlinear mapping $S$ in the following way: for any
$(u,v)\in X_{T,R}$, $(\bar{u},\bar{v})=S(u,v)$ is the unique
solution of the following linear wave equation with nonnegative
potential
\begin{eqnarray}
\left\{
\begin{array}{l}
\alpha(\bar{u}_{tt}-\Delta
\bar{u}+m_1^2\bar{u}+K_1(x)\bar{u})=a_2^\prime(p+1)|v|^{q+1}|u|^{p-1}u\\
\bar{v}_{tt}-\Delta
\bar{v}+m_2^2\bar{v}+K_2(x)\bar{v}=a_2^\prime(q+1)|u|^{p+1}|v|^{q-1}v\\
\bar{u}(0,x)=u_0, \bar{v}(0,x)=v_0\\
\bar{u}_t(0,x)=u_1, \bar{v}_t(0,x)=v_1
\end{array}
\right.\label{CKGE_P3}
\end{eqnarray}

Obviously, by Theorem 2.2 the existence and uniqueness of the
solution $(\bar{u}(t),\bar{v}(t))$ can be obtained for $(u,v)\in
H^1_{0, K_1}\times H^1_{0, K_2}$.

We next claim that, for suitable $R$ and $T$, $S$ is a contraction
mapping satisfying $S(X_{T,R})\subseteq X_{T,R}$. Indeed, in the
following part of the proof, we will take $R^2=e(u_0,v_0)$.
 Given $(u,v)\in
X_{T,R}$, for every $t\in(0, T]$ the corresponding solution
$(\bar{u},\bar{v})=S(u,v)$ satisfies the energy identity:
\begin{eqnarray*}
&&\frac{\alpha}{2}\left(\|\bar u_t(t)\|^2+\|\nabla \bar
u(t)\|^2+m_1^2\|\bar{u}(t)\|^2+\int
K_1(x)|\bar{u}(t)|^2dx\right)\\
&&\hskip20pt=\frac{\alpha}{2}(\|u_0\|_{H^1_{0,
K_1}}^2+\|u_1\|_2^2)
+a_2^\prime(p+1)\int _0^t\int|v|^{q+1}|u|^{p-1}u\bar{u}d\tau,\\
&&\frac{1}{2}\left(\|\bar v_t(t)\|^2+\|\nabla \bar
v(t)\|^2+m_2^2\|\bar{v}(t)\|^2+\int
K_2(x)|\bar{v}(t)|^2dx\right)\\
&&\hskip20pt=\frac{1}{2}(\|v_0\|_{H^1_{0,
K_2}}^2+\|v_1\|_2^2)+a_2^\prime(q+1)\int _0^t\int|u|^{p+1}|v|^{q-1}v\bar{v}d\tau.
\end{eqnarray*}
For the last terms of the right hand side above, we have by Lemma
2.1
\begin{eqnarray*}
\int _0^T\int|v|^{q+1}|u|^{p-1}u\bar{u}d\tau\leq cT
R^{2(p+q+1)}+2\int_0^T\|\nabla \bar{u}\|^2d\tau,\\
\int _0^T\int|u|^{p+1}|v|^{q-1}v\bar{v}d\tau\leq cT
R^{2(p+q+1)}+2\int_0^T\|\nabla \bar{u}\|^2d\tau.
\end{eqnarray*}
Thus taking the maximum on $[0, T]$, we have
\begin{eqnarray*}
e(\bar u,\bar v)\leq \frac{1}{2}R^2+cTR^{2(p+q+1)}.
\end{eqnarray*}
Obviously, taking $T>0$ sufficient small, we have $e(\bar u,\bar
v)\leq R^2$, which means $S(X_{T,R})\subseteq X_{T,R}$.

We next show that $S$ is a contraction in $X_{T,R}$ with the
distance $e(u_1-u_2, v_1-v_2)$. Take $(u_1,v_1)$ and $(u_2,v_2)$
from $X_{T,R}$, and denote the corresponding solution of
(\ref{CKGE_P3}) by $(\bar{u}_1,\bar{v}_1)$ and
$(\bar{u}_2,\bar{v}_2)$, respectively. Then by mean value theorem we
have
\begin{eqnarray*}
\left||v_1|^{q+1}|u_1|^{p-1}u_1-|v_2|^{q+1}|u_2|^{p-1}u_2\right|\leq
|v_1-v_2|\omega_1(x,t)|u_1|^p+|u_1-u_2|\omega_2(t,x)|v_2|^{q+1},\\
\left||u_1|^{p+1}|v_1|^{q-1}v_1-|u_2|^{p+1}|v_2|^{q-1}v_2\right|\leq
|u_1-u_2|\omega_3(x,t)|v_1|^q+|v_1-v_2|\omega_4(t,x)|u_2|^{p+1}.
\end{eqnarray*}
where
\begin{eqnarray*}
&&\omega_1(t)\leq(|v_1|+|v_2|)^q,\\
&&\omega_2(t)\leq(|u_1|+|u_2|)^{p-1},\\
&&\omega_3(t)\leq(|u_1|+|u_2|)^p,\\
&&\omega_4(t)\leq(|v_1|+|v_2|)^{q-1}.
\end{eqnarray*}

Thus by H\"older inequality and Sobolev inequality we have
\begin{eqnarray*}
e(\bar{u}_1-\bar{u}_2,\bar{v}_1-\bar{v}_2)\leq cR^{2(p+q)}T
e(u_1-u_2,v_2-v_2).
\end{eqnarray*}
If we let $T$ sufficiently small, then $cR^{2(p+q)}T<1$. This
implies that $S$ is a contraction in $X_{T,R}$.

 From the above argument, by applying Banach fixed point theorem, we
obtain the existence of a unique solution of the system
(\ref{CKGE_P1}) on a maximum time interval.

For the last statement, let the interval $[0, T_{\max})$ is maximal
interval where the solution of (\ref{CKGE_P1}) exists. We assume
that $T_{\max}<\infty$ and $K=\lim_{t\rightarrow
T_{\max}^-}(\alpha\|u(t)\|_2^2+\|v\|^2_2)<\infty$. Then there exists
a sequence $\{t_j\}_{j=1}^\infty$ such that
\begin{eqnarray*}
&&t_j\rightarrow T_{\max}^-,\\
&&\alpha\|u(t_j)\|_2^2+\|v(t_j)\|^2_2\leq K.
\end{eqnarray*}

Using the same argument as above with the initial data at $t_j$ we see
that there exists a unique solution of (\ref{CKGE_P1}) on $[t_j,
t_j+T_{\max,j}^-)$. Thus we can get $T_{\max}<t_j+T_{\max,j}$ for
some $j$ large enough. Obviously this contradicts the definition
of $T_{\max}$.

Thus we have completed the proof of the local existence theorem.
\begin{flushright}
$\Box$
\end{flushright}

\section{Standing wave with ground state}

If a real function $(\phi, \psi)$ verifies the following system
\begin{eqnarray}
\left\{\begin{array}{l}
-\alpha\Delta\phi+\alpha m_1^2\phi+\alpha K_1(x)\phi=a_2^\prime(p+1)|\psi|^{q+1}|\phi|^{p-1}\phi\\
-\Delta\psi+m_2^2\psi+K_2(x)\psi=a_2^\prime(q+1)|\phi|^{p+1}|\psi|^{q-1}\psi
\end{array}
\right.\label{Euler_eq}
\end{eqnarray}
and
\begin{eqnarray*}
(\phi,\psi)\in H^1_{0, K_1}(\R^n)\times H^1_{0,K_2}(\R^n),
\end{eqnarray*}
then $(u,v)=(\phi,\psi)$ verifies the system (\ref{CKGE_P1}) for
$t\geq 0, x\in\R^n$.

We now define the action $J(\phi,\psi)$ of the solution $(\phi,
\psi)$ of the system (\ref{Euler_eq}) as follows
\begin{eqnarray}
J(\phi,\psi)&=&\frac{1}{2}\int\left(\alpha|\nabla \phi(x)|^2+|\nabla
\psi(x)|^2+\alpha m_1^2|\phi(x)|^2+m_2^2|\psi(x)|^2\right.\nonumber\\
&&+\alpha
\left.K_1(x)|\phi(x)|^2+K_2(x)|\psi(x)|^2\right)dx-a_2^\prime\int|\phi(x)|^{p+1}|\psi(x)|^{q+1}dx.
\label{J_def}
\end{eqnarray}
 In
addition, we let
\begin{eqnarray}
I(\phi,\psi)&=&\int\left(\alpha|\nabla \phi(x)|^2+|\nabla
\psi(x)|^2+\alpha m_1^2|\phi(x)|^2+m_2^2|\psi(x)|^2\right.\nonumber\\
&&+\left.\alpha
K_1(x)|\phi(x)|^2+K_2(x)|\psi(x)|^2\right)dx-a_2^\prime(p+q+2)\int|\phi(x)|^{p+1}|\psi(x)|^{q+1}dx.\label{I_def}
\end{eqnarray}

We now state a proposition, which describes the relation between
$J(\phi,\psi)$ and $I(\phi,\psi)$.
\begin{proposition}
Let $(\phi,\psi)\in H^1_{0, K_1}\times H^1_{0, K_2}$ satisfy $(\phi,\psi)\neq0$ and set
$(\phi^\lambda(x),\psi^\lambda(x))=(\lambda\phi(x),\lambda\psi(x))$ for $\lambda>0$.
Then for $p>0$ and $q>0$ there exists a unique $\lambda_1$ such
that
\begin{eqnarray}
I(\phi^{\lambda},\psi^{\lambda})\left\{
\begin{array}{l}
>0 \mathrm{\ when \ } \lambda\in(0,\lambda_1),\\
=0 \mathrm{\ when\ }\lambda=\lambda_1, \\
<0 \mathrm{\ when\ }\lambda\in(\lambda_1,\infty).
\end{array}
\right.\label{I_1}
\end{eqnarray}
and
\begin{eqnarray}
J(\phi^\lambda,\psi^\lambda)\leq
J(\phi^{\lambda_1},\psi^{\lambda_1})\label{J-2}
\end{eqnarray}
for all $\lambda>0$.
\end{proposition}
\textbf{Proof.} By (\ref{I_def}) and
$(\phi^\lambda,\psi^\lambda)=(\lambda\phi,\lambda\psi)$, we easily
see that $I(\phi^\lambda,\psi^\lambda)$ is continuous in $\lambda$.
Thus it is obvious that there exists a unique $\lambda_1$ such that
$I(\phi^\lambda,\psi^\lambda)$ satisfies the relation (\ref{I_1}).

Furthermore, by a direct computation we see that for $\lambda>0$
\begin{eqnarray*}
\lambda\frac{d}{d\lambda}J(\phi^\lambda,\psi^\lambda)=I(\phi^\lambda,\psi^\lambda),
\end{eqnarray*}
which implies the property (\ref{J-2}) by (\ref{I_1}).
\begin{flushright}
$\Box$
\end{flushright}

 Next define the set $M$ by
\begin{eqnarray}
M=\{(\phi,\psi)\in H_1(\R^n)\times H_2(\R^n); I(\phi,\psi)=0, (\phi,
\psi)\ne 0\}.\label{M_def}
\end{eqnarray}
We then set the following constrained variational problem
\begin{eqnarray*}
d=\min_{(\phi,\psi)\in M} J(\phi,\psi).
\end{eqnarray*}

Now we are in a position to state the theorem about the ground state
of (\ref{Euler_eq}).
\begin{theorem} Assume that $1<p,q<\infty$ when $n=2$ and $\displaystyle 1<p,q<\frac{2}{n-2}$ when $n=3$.
 Then there exists $(\Phi,\Psi)\in M$ such that

(1) $\displaystyle J(\Phi, \Psi)=\inf_{(\phi, \psi)\in M}J(\phi,
\psi)=d$;

(2) $(\Phi, \Psi)$ is a ground state solution of (\ref{Euler_eq}).
\end{theorem}

To prove the above theorem we first introduce a compactness lemma,
whose proof can be found in \cite{23Omana}, \cite{22Zhang}.
\begin{lemma}
Let $\displaystyle 1\leq r<\frac{n+2}{n-2}$ when $n\geq 3$ and
$1\leq r<\infty$ when $n=1, 2$. Then for $K(x)\geq 0$ ($\forall
x\in\R^n$) the embedding $H^1_{0, K}\hookrightarrow L^{r+1}$ is
compact.
\end{lemma}

\textbf{Proof of Theorem 3.1.} Firstly we claim that $J(\phi, \psi)$
is bounded below on $M$. Indeed, By (\ref{J_def}) and (\ref{M_def})
we see that
\begin{eqnarray}
J(\phi,\psi)&=&\frac{p+q}{2(p+q+2)}\int \left(\alpha|\nabla
\phi|^2+|\nabla
\psi| +\alpha m_1^2|\phi|^2+m_2^2|\psi|^2+\alpha K_1|\phi|^2+K_2|\psi|^2\right)dx\nonumber\\
&\geq& 0.\label{J-BelowB}
\end{eqnarray}

Then there exists a minimizing sequence $\{(\phi_j,\psi_j)\}_{j=1}^\infty$ satisfying
\begin{eqnarray}
&&I(\phi_j,\psi_j)=0 \mathrm{\ for\ every\ }j\in\N,\label{I_assm}\\
&&J(\phi_j,\psi_j)\rightarrow d \mathrm{\ as\
}j\rightarrow\infty.\label{J_assm}
\end{eqnarray}
By (\ref{J-BelowB}) and (\ref{J_assm}) there exists a
subsequence of $\{(\phi_j,\psi_j)\}_{j=1}^\infty$, which we still
denote by $\{(\phi_j,\psi_j)\}_{j=1}^\infty$, such that
\begin{eqnarray*}
(\phi_j,\psi_j)\rightharpoonup (\phi_\infty,\psi_\infty) \mathrm{\
weakly\ in \ } H_{0, K_1}^1\times H_{0, K_2}^1.
\end{eqnarray*}

By Lemma 3.1 we have
\begin{eqnarray}
&&(\phi_j,\psi_j)\rightarrow (\phi_\infty,\psi_\infty) \mathrm{\
strongly\ in\ } L^2(\R^n)\times L^2(\R^n);\label{TH21-1}\\
&&(\phi_j,\psi_j)\rightarrow (\phi_\infty,\psi_\infty) \mathrm{\
strongly\ in\ } L^{2(p+1)}(\R^n)\times
L^{2(q+1)}(\R^n).\label{TH21-2}
\end{eqnarray}


Next by a contradiction argument we prove that
\begin{eqnarray}
(\phi_\infty,\psi_\infty)\not\equiv 0.\label{TH21-3}
\end{eqnarray}

We assume that $(\phi_\infty,\psi_\infty)\equiv 0$. Then by
(\ref{TH21-1}) and (\ref{TH21-2}) we obtain
\begin{eqnarray}
&&(\phi_j,\psi_j)\rightarrow 0 \mathrm{\ in\ } L^2(\R^2)\times L^2(\R^2);\label{TH21-4}\\
&&(\phi_j,\psi_j)\rightarrow 0 \mathrm{\ in\ }
L^{2(p+1)}(\R^2)\times L^{2(q+1)}(\R^2).\label{TH21-5}
\end{eqnarray}
Noting the fact $(\phi_j,\psi_j)\in M$, we see that
$I(\phi_j,\psi_j)=0$ implies that
\begin{eqnarray*}
\int (\alpha(|\nabla \phi_j|^2+m_1^2|\phi|^2+K_1(x)|\phi|^2)+|\nabla
\psi|^2+m_2^2|\psi|^2+K_2(x)|\psi|^2)dx\rightarrow 0,
\end{eqnarray*}
which means
\begin{eqnarray}
\int (\alpha(|\nabla \phi_j|^2+K_1(x)|\phi_j|^2)+|\nabla
\psi_j|^2+K_2(x)|\psi_j|^2)dx\rightarrow
0\label{Th21_Contradiction1}
\end{eqnarray}
as $j\rightarrow \infty$.

On the other hand, from (\ref{J_def}) and (\ref{J_assm}) we see
\begin{eqnarray}
\int (\alpha(|\nabla \phi_j|+K_1(x)|\phi_j|^2)+|\nabla
\psi_j|+K_2(x)|\psi_j|^2)\rightarrow
\frac{2(p+q+2)}{p+q}d\label{Th21_Contradiction2}
\end{eqnarray}
as $j\to\infty$.

Obviously if we can show $d>0$, then there exists a contradiction
between (\ref{Th21_Contradiction1}) and (\ref{Th21_Contradiction2}),
which implies that it is impossible that
$(\phi_\infty,\psi_\infty)\equiv 0$. We next show it.

Indeed, by (\ref{J_def}), (\ref{I_def}) and (\ref{M_def}) we have
\begin{eqnarray}
&&\hskip-34ptJ(\phi_j,\psi_j)\nonumber\\
&&\hskip-28pt=\frac{p+q}{2(p+q+2)}\int (\alpha(|\nabla
\phi_j|+m_1^2|\phi_j|^2+K_1(x)|\phi_j|^2)+|\nabla
\psi_j|+m_2^2|\psi_j|^2+K_2(x)|\psi_j|^2)dx.\label{J-1}
\end{eqnarray}

Using H\"older inequality, the Sobolev`s embedding theorem and
$I(\phi_j,\psi_j)=0$ we obtain
\begin{eqnarray*}
\int |\phi_j|^{p+1}|\psi_j|^{q+1}dx&\leq& \left(\int
|\phi_j|^{2(p+1)}dx\right)^{1/2}
\left(\int |\psi_j|^{2(q+1)}dx\right)^{1/2}\nonumber\\
&\leq& c\left(\alpha\int(|\nabla
\phi_j|^2+m_1^2|\phi_j|^2)dx\right)^{(p+1)/2}\left(\int(|\nabla
\psi_j|^2+m_2^2|\psi_j|^2)dx\right)^{(q+1)/2}
\nonumber\\
&\leq& c\left(\int(\alpha(|\nabla \phi_j|^2+m_1^2|\phi_j|^2)+|\nabla
\psi_j|^2+m_2^2|\psi_j|^2)dx\right)^{(p+q+2)/2}
\end{eqnarray*}
for some constant $c>0$, and
\begin{eqnarray*}
&&a_2^\prime(p+q+2)\int|\phi_j|^{p+1}|\psi_j|^{q+1}dx\\
&&=\int (\alpha(|\nabla
\phi_j|^2+m_1^2|\phi_j|^2+K_1(x)|\phi_j|^2)+|\nabla
\psi_j|^2+m_2^2|\psi_j|^2+K_2(x)|\psi_j|^2)dx
\end{eqnarray*}
for every $j\in\N$.

Thus we see that
\begin{eqnarray*}
&&\int (\alpha(|\nabla
\phi_j|^2+m_1^2|\phi_j|^2+K_1(x)|\phi_j|^2)+|\nabla
\psi_j|^2+m_2^2|\psi_j|^2+K_2(x)|\psi_j|^2)dx\\
&&\leq c'\left(\int(\alpha(|\nabla \phi_j|^2+m_1^2|\phi_j|^2)+|\nabla
\psi_j|^2+m_2^2|\psi_j|^2)dx\right)^{(p+q+2)/2}\\
&&\leq c'\left(\int(\alpha(|\nabla
\phi_j|^2+m_1^2|\phi_j|^2+K_1|\phi_j|^2)+|\nabla
\psi_j|^2+m_2^2|\psi_j|^2+K_2|\psi_j|^2)dx\right)^{(p+q+2)/2},
\end{eqnarray*}
for some constant $c'>0$, which implies that for some constant $c''>0$
\begin{eqnarray*}
\int(\alpha(|\nabla \phi_j|^2+m_1^2|\phi_j|^2+K_1|\phi_j|^2)+|\nabla
\psi_j|^2+m_2^2|\psi_j|^2+K_2|\psi_j|^2)dx\geq c''>0
\end{eqnarray*}
for every $j\in\N$.

Therefore from (\ref{J-1}) it follows that
\begin{eqnarray*}
d=\lim_{j\rightarrow\infty}J(\phi_j,\psi_j)>0.
\end{eqnarray*}

Thus we have completed the proof of
$(\phi_\infty,\psi_\infty)\not\equiv 0$.

By Proposition 3.1 there exists a unique $\lambda_0\in(0,\infty)$
such that $I(\phi_\infty^{\lambda_0},\psi_\infty^{\lambda_0})=0$,
where
$$
(\phi^\lambda(x),\psi^\lambda(x))=(\lambda\phi(x),\lambda\psi(x)).
$$
Then it follows from Proposition 3.1 and $I(\phi_j,\psi_j)=0$ that
\begin{eqnarray}
J(\phi_\infty^{\lambda_0},\psi_\infty^{\lambda_0})\leq J(\phi_j,\psi_j)
\rightarrow d\nonumber
\end{eqnarray}
as $j\rightarrow\infty$.

Noting the fact
$I(\phi_\infty^{\lambda_0},\psi_\infty^{\lambda_0})=0$,
we have $I(\Phi,\Psi)=0$ and $J(\Phi,\Psi)=d$ with
$\Phi=\phi_\infty^{\lambda_0}$ and $\Psi=\psi_\infty^{\lambda_0}$.

Since $(\Phi,\Psi)$ is a solution of Lagrange-Euler equation
(\ref{Euler_eq}), there exists a Lagrange multiplier $\Theta$ such
that
\begin{eqnarray}
J_\Phi(\Phi,\Psi)+\Theta I_\Phi(\Phi,\Psi)=0,\label{S3-1}\\
J_\Psi(\Phi,\Psi)+\Theta I_\Psi(\Phi,\Psi)=0,\label{S3-2}
\end{eqnarray}
and
\begin{eqnarray}
\langle J_\Phi(\Phi,\Psi)+\Theta I_\Phi(\Phi,\Psi),\Phi\rangle=0,\label{S3-11}\\
\langle J_\Psi(\Phi,\Psi)+\Theta I_\Psi(\Phi,\Psi),\Psi\rangle=0.\label{S3-21}
\end{eqnarray}
Since $I(\Phi,\Psi)=0$, it follows from (\ref{S3-1})-(\ref{S3-21})
that
\begin{eqnarray*}
\Theta\int |\Phi|^{q+1}|\Psi|^{p+1}dx=0
\end{eqnarray*}
which implies $\Theta\equiv 0$. Therefore $(\Phi,\Psi)$ solves the
equation (\ref{Euler_eq}) in $H^1_{0,K_1}(\R^n)\times H^1_{0,
K_2}(\R^n)$. This completes the proof of Theorem 3.1.\hfill$\Box$

\begin{remark}
The standing wave of the system (\ref{CKGE_P1}) is the form of
($\exp(i\omega_1 t)\phi_{\omega_1}(x) ,
\exp(i\omega_2t)\psi_{\omega_2}(x)$), where $(\phi_{\omega_1},
\psi_{\omega_2})$ is the ground state solution of the corresponding
Lagrange-Euler equations. Indeed, in the above argument we just
consider the case $(\omega_1,\omega_2)=0$. Recently, for the single
Klein-Gorodn equation (\ref{SKGE}) with $K(x)\equiv 0$, the strong
instability of the standing waves $\exp(i\omega
t)\varphi_{\omega}(x)$, where $\varphi_{\omega}(x)$ is a ground
state solution for the corresponding single Lagrange-Euler equation,
has been considered in \cite{26Liu}, \cite{13Ohta} with $|\omega|<1$.
But for $(\omega_1,\omega_2)\ne 0$, the strong instability of the
standing waves of the system (\ref{CKGE_P}) is unknown.
\end{remark}

\section{Blow up with low initial energy when $n=2,3$}
In the section, based on the result obtained in Section 3 we prove
the blow up result by potential well argument and concavity method,
which leads to the instability of standing waves of the system
(\ref{CKGE_P1}).

Firstly define two sets $\Gamma_1$ and $\Gamma_2$ as
\begin{eqnarray*}
\Gamma_1=\{(u,v)\in H_1\times H_2; J(u,v)<d, I(u,v)<0\},\\
\Gamma_2=\{(u,v)\in H_1\times H_2; J(u,v)<d, I(u,v)>0\},
\end{eqnarray*}
where $d$ is defined in theorem 3.1.

\begin{lemma}
Assume $1<p,q<\infty$ when $n=2$ and $\displaystyle
1<p,q<\frac{2}{n-2}$ when $n=3$. Let the initial energy satisfying
$E(0)<d$. Then the set $\Gamma_1$ is invariant under the flow
generated by (\ref{CKGE_P1}) in the sense that: If nonzero $(u_0,
v_0)\in \Gamma_1$, then the unique solution $(u(t),v(t))$ of the
equations (\ref{CKGE_P1}), $0\leq t<T_{\max}$, with the initial data
$(u_0, v_0)$ satisfies
\begin{eqnarray*}
(u(t),v(t))\in \Gamma_1, \mathrm{\ for\ }t\in[0, T_{\max}),
\end{eqnarray*}
where $T_{\max}>0$ is the maximum existing time of the solution
$(u(t),v(t))$.
\end{lemma}
\textbf{Proof.} It is observed by the conservation law of energy
(\ref{E_conservation}) that
\begin{eqnarray*}
J(u(t),v(t))\leq E(t)=E(0)<d
\end{eqnarray*}
for $0\leq t<T_{\max}$.

To prove $(u(t),v(t))\in\Gamma_1$, there is only one thing left to be
checked: $I(u(t),v(t))<0$ for $0\leq t<T_{\max}$. In the following we
show it by a contradiction argument.

We assume that it is wrong that $I(u(t),v(t))<0$ for $0\leq t<T_{\max}$. Then by
continuity, we see there exists a time $T>0$ such that
\begin{eqnarray*}
T=\min\{0<t<T_{\max}; I(u(t),v(t))=0\}.
\end{eqnarray*}
It is natural that $I(u(T), v(T))=0$ by the continuity of $I(u(t),
v(t))$ in $t$. Then $(u(T),v(T))\in M$. By Theorem 3.1 we see that
it is impossible that $J(u(T),v(T))<d$ and $(u(T),v(T))\in M$. Thus,
we have obtained that $I(u(t),v(t))<0$ for $0\leq t<T_{\max}$. So
$\Gamma_1$ is invariant under the flow generated by (\ref{CKGE_P1}).
\begin{flushright}
$\Box$
\end{flushright}
\begin{theorem}
Assume that $1<p,q<\infty$ when $n=2$ and $\displaystyle
1<p,q<\frac{2}{n-2}$ when $n=3$. If the initial datum $(u_0,v_0)$
and $(u_1,v_1)$ satisfy $E(0)<d$ and $I(u_0, v_0)<0$, then the
solution $(u(t),v(t))$ of the Cauchy problem (\ref{CKGE_P1}) blows
up in a finite time, that is,
\begin{eqnarray*}
\lim_{t\rightarrow T_{\max}^-}(\|u(t)\|^2+\|v(t)\|^2)=\infty.
\end{eqnarray*}
\end{theorem}
\textbf{Proof.} By Lemma 4.1 we see that
\begin{eqnarray}
I(u(t),v(t))<0\label{I<0}
\end{eqnarray}
for every $t\in[0, T_{\max})$.

We first define the following auxiliary function
\begin{eqnarray}
G(t)=\int (\alpha|u(t,x)|^2+|v(t,x)|^2)dx+b(t+T_1)^2.\label{G_def}
\end{eqnarray}
where $b$ and $T_1$ is two positive parameters, which will be
determined later.

By simple calculation we see that
\begin{eqnarray}
G^\prime(t)&=&\frac{d}{dt}G(t)\nonumber\\
&=&2\int (\alpha
u(t,x)u_t(t,x)+v(t,x)v_t(t,x))dx+2b(t+T_1),\label{G1diff}
\end{eqnarray}
and
\begin{eqnarray}
\frac{1}{2}G^{\prime\prime}(t)&=&\int
(\alpha|u_t(t,x)|^2+|v_t(t,x)|^2)dx+\int \alpha
u(t,x)u_{tt}(t,x)dx+\int
v(t,x)v_{tt}(t,x)dx+b\nonumber\\
&=&\int (\alpha|u_t(t,x)|^2+|v_t(t,x)|^2)dx+\int \alpha
u(t,x)(\Delta
u-m_1^2u-K_1(x)u+a_1|v|^{q+1}|u|^{p-1}u)dx\nonumber\\
&&+\int v(t,x)(\Delta
v-m_2^2v-K_2(x)v+a_2|u|^{p+1}|v|^{q-1}v)dx+b\nonumber\\
&=&\int (\alpha|u_t(t,x)|^2+|v_t(t,x)|^2)dx+a_2^\prime(p+q+2)\int
|v|^{q+1}|u|^{p+1}dx\nonumber\\
&& -\alpha\int(|\nabla u|^2+m_1^2|u|^2+K_1(x)|u|^2)dx -\int(|\nabla
v|^2+m_2^2|v|^2+K_2(x)|v|^2)dx+b\label{G2diff}\\
&=&\int (\alpha|u_t(t,x)|^2+|v_t(t,x)|^2)dx-I(u(t),v(t))+b\nonumber.
\end{eqnarray}

Thus it is obvious by (\ref{I<0}) that
\begin{eqnarray}
G^{\prime\prime}(t)>0\label{G2diff_p}
\end{eqnarray}
on $[0, T_{\max})$, which implies that the function $G(t)$ is
convex in $t$.

 From (\ref{E_conservation}), (\ref{E_def}) and (\ref{G2diff}) it
follows that
\begin{eqnarray*}
G^{\prime\prime}(t)&=&(p+q+4)\int
(\alpha|u_t|^2+|v_t|^2)dx+(p+q)\left[\alpha\int(|\nabla
u|^2+m_1^2|u|^2+K_1(x)|u|^2)dx\right. \nonumber\\
&&+\left.\int(|\nabla
v|^2+m_1^2|v|^2+K_2(x)|v|^2)dx\right]-2(p+q+2)E(t)+2b\nonumber\\
&\geq&(p+q+4)\int
(\alpha|u_t|^2+|v_t|^2)dx+(p+q)\min\{m_1,m_2\}G(t)-2(p+q+2)E(t)+2b.
\end{eqnarray*}

We now choose a sufficiently large time $T_1$ and suitable $b$
satisfying
\begin{eqnarray}
G^\prime(0)>0\label{G1diff(0)}
\end{eqnarray}
and 
\begin{eqnarray}
(p+q)\min\{m_1,m_2\}G(0)-2(p+q+2)E(0)>(p+q+2)b.
\end{eqnarray}

Thus we see that $G(t)$ is strictly increasing on $[0, T_{\max})$ by
(\ref{G2diff_p}) and (\ref{G1diff(0)}). So we have
\begin{eqnarray*}
(p+q)\min\{m_1,m_2\}G(t)-2(p+q+2)E(0)> (p+q+2)b
\end{eqnarray*}
for every $t\in [0, T_{\max})$.

Then for every $t\in[0, T_{\max})$ we get
\begin{eqnarray}
G^{\prime\prime}(t)\geq (p+q+4)\left[\int
(\alpha|u_t|^2+|v_t|^2)dx+b\right].\label{G2diff>}
\end{eqnarray}

By (\ref{G_def}), (\ref{G1diff}) and (\ref{G2diff>}) we easily see
that
\begin{eqnarray}
G^{\prime\prime}(t)G(t)-\frac{p+q+4}{4}(G^\prime(t))^2&\geq&
(p+q+2)\left\{\left[\int (\alpha
|u_t|^2+|v_t|^2)dx+b\right]\right.\nonumber\\
&&\times\left[\int(\alpha|u|^2+|v|^2)dx+b(t+T_1)^2\right]\nonumber\\
&&-\left.\left[\int(\alpha
u_tu+v_tv)dx+b(t+T_1)\right]^2\right\}\nonumber\\
&\geq& 0,\label{G2ttG-G1t>0}
\end{eqnarray}
where the last inequality comes from Cauchy-Schwartz inequality.

 And by direct computation we have
 \begin{eqnarray}
 \frac{d}{dt}G^{-(p+q)/4}(t)=-\frac{p+q}{4}G^{-(p+q+4)/4}(t)G^\prime(t)<0\label{G_con01}
 \end{eqnarray}
 and
\begin{eqnarray}
\frac{d^2}{dt^2}G^{-(p+q)/4}(t)&=& -\frac{p+q}{4}G^{-(p+q+8)/4}
\left(G^{\prime\prime}(t)G(t)-\frac{p+q+4}{4}(G^\prime(t))\right)\nonumber\\
&\leq& 0.\label{G_con0}
\end{eqnarray}
for every $t\in[0, T_{\max})$, which implies that $G^{-(p+q)/4}(t)$
is concave for $t\geq 0$. From (\ref{G_con01}) and (\ref{G_con0}) it follows that
the function $G^{-(p+q)/4}\rightarrow 0$ when $t<T_{\max}$ and
$t\rightarrow T_{\max}$ ($\displaystyle
T_{\max}\leq\frac{4G(0)}{(p+q)G^\prime(0)}$). Thus we see that there exists
a finite time $T_{\max}>0$ such that
\begin{eqnarray*}
\lim_{t\rightarrow T^-_{\max}}(\alpha\|u(t)\|_2^2+\|v(t)\|_2^2)=\infty.
\end{eqnarray*}
\begin{flushright}
$\Box$
\end{flushright}
Next we will consider the instability of standing waves of the
equation (\ref{CKGE_P1}). Before we give the instability
theorem, we introduce the definition of the instability of standing
waves.
\begin{definition}
The standing wave $(u(t,x),v(t,x))=(\phi(x),\psi(x))$ of the system
(\ref{CKGE_P1}) is instable by blow up if for any $\epsilon>0$ there
exists $(u_0, v_0)\in H^1_{0,K_1}\times H^1_{0,K_2}$ such that
$\|(u_0,v_0)-(\phi,\psi)\|_{H^1_{0,K_1}\times H^1_{0,K_2}}<\delta$
and the corresponding solution of the system (\ref{CKGE_P}) blows up
in a finite time with the initial data
\begin{eqnarray*}
\left\{\begin{array}{l}
u(0,x)=u_0, u_t(0,x)=0,\\
v(0,x)=v_0, v_t(0,x)=0.
\end{array}
\right.
\end{eqnarray*}
\end{definition}

Now we are in a position to state the instability theorem.
\begin{theorem}Assume that $1<p, q<\infty$ when $n=2$ and $\displaystyle 1<p,
q<\frac{2}{n-2}$ when $n=3$. Let $(\phi,\psi)$ be a ground state
solution of (\ref{Euler_eq}). Then for any $\epsilon>0$ there exists
$T_{\max}<\infty$ and $(u_0,v_0)\in H^1_1\times H^1_2$ with the
property
\begin{eqnarray*}
\|u_0-\phi\|_{H^1_{0, K_1}}<\epsilon,\\
\|v_0-\psi\|_{H^1_{0, K_2}}<\epsilon.
\end{eqnarray*}
such that the solution $(u,v)$ of the system (\ref{CKGE_P}) blows up
in a finite time $T_{\max}$ with the initial data
\begin{eqnarray}
\left\{\begin{array}{l}
u(0,x)=u_0, u_t(0,x)=0\\
v(0,x)=v_0, v_t(0,x)=0
\end{array}
\right.\label{Initi-data}
\end{eqnarray}
\end{theorem}
\textbf{Proof.} By (\ref{E_def}) and (\ref{Initi-data}) we easily
see that
\begin{eqnarray}
E(0)=J(u_0, v_0).\label{TH32-1}
\end{eqnarray}
We next let
\begin{eqnarray}
u_0(x)=\gamma \phi(x), v_0(x)=\gamma\psi(x)\label{Initi-data1}
\end{eqnarray}
where $\gamma>1$.

Obviously, for any $\epsilon>0$ we can take a suitable $\gamma>1$
such that
\begin{eqnarray*}
\|u_0-\phi\|_{H^1_{0,K_1}}=(\gamma-1)\|\phi\|_{H^1_{0,K_1}}<\epsilon,\\
\|v_0-\psi\|_{H^1_{0,K_2}}=(\gamma-1)\|\psi\|_{H^1_{0,K_2}}<\epsilon.
\end{eqnarray*}
Noting the fact $\gamma>1$, we see by (\ref{Initi-data1}) and Proposition
3.1 that
\begin{eqnarray*}
I(u_0,v_0)<I(\phi,\psi)=0\\
J(u_0,v_0)<J(\phi,\psi)=d
\end{eqnarray*}
and by (\ref{TH32-1}) we have $E(0)<d$.

Thus by Theorem 3.1 we have completed the proof of Theorem 3.2.
\begin{flushright}
$\Box$
\end{flushright}
As a byproduct we have the following global existence theorem for
the system (\ref{CKGE_P}).
\begin{theorem} Assume that $1<p, q<\infty$ when $n=2$ and $\displaystyle 1<p,
q<\frac{2}{n-2}$ when $n=3$. If the initial data satisfy that
$E(0)<d$ and $I(u_0, v_0)>0$ then the corresponding solution
$(u(t,x),v(t,x))$ of the system (\ref{CKGE_P1}) exists globally.
\end{theorem}
\textbf{Proof.} As in the proof of Lemma 3.1 we can prove that
$\Gamma_2$ is invariant under the flow generated by the system
(\ref{CKGE_P1}). Thus we see that
\begin{eqnarray}
J(u(t),v(t))<d,\label{TH33-2}\\
I(u(t),v(t))>0\label{TH33-3}
\end{eqnarray}
for each $t\in[0, T_{\max})$.

Thus by (\ref{J_def}) we get
\begin{eqnarray}
J(u(t),v(t))&>&\frac{p+q}{2(p+q+2)}\left[\int \alpha(|\nabla
u(t)|^2+m_1^2|u(t)|+K_1(x)|u(t)|^2)dx\right.\nonumber\\
&&+\left.\int (|\nabla
v(t)|^2+m_2^2|v(t)|+K_2(x)|v(t)|^2)dx\right]\nonumber\\
&\geq& 0.\label{TH33-1}
\end{eqnarray}
 From (\ref{TH33-2}) and (\ref{TH33-1}) it follows that
\begin{eqnarray*}
\int (\alpha|u_t(t)|^2+|v_t(t)|^2)dx<\frac{2(p+q+2)}{p+q}d
\end{eqnarray*}
for each $t\in[0, T_{\max})$.

By (\ref{TH33-2}), (\ref{TH33-3}) and (\ref{TH33-1}) we see that
\begin{eqnarray*}
&&\int (\alpha|\nabla u|^2+|\nabla v|)dx<\frac{2(p+q+2)}{p+q}d,\\
&&\int (\alpha m_1^2|u|^2+m_2^2|v|^2)dx<\frac{2(p+q+2)}{p+q}d,\\
&&\int(\alpha K_1(x)|u|^2+K_2|v|^2)dx<\frac{2(p+q+2)}{p+q}d
\end{eqnarray*}
for each $t\in [0, T_{\max})$.

Thus we see that the solution is uniformly bounded on $[0, T_{\max})$. The proof of the theorem is finished.
\begin{flushright}
$\Box$
\end{flushright}

\section{Blow up with non-positive initial energy}
Because of the embedding theorem (Lemma 3.1) we cannot claim the
blowing up result as in Section 4 for other case. So in this section
we will prove the blow up result for the system (\ref{CKGE_P1}) by
the concavity method when the initial energy is not positive. We
state our theorem.
\begin{theorem} Assume that $p$ and $q$ satisfy the conditions (\ref{pq-cond1}), (\ref{pq-cond2}) and (\ref{pq-cond3}).
If the nonzero datum $(u_0,v_0)\in H^1_{0, K_1}\times H^1_{0, K_2}$
and $(u_1,v_1)\in L^2\times L^2$ satisfy
\begin{eqnarray*}
E(0)<0,
\end{eqnarray*}
or
\begin{eqnarray*}
\int \alpha u_0u_1+v_0v_1\geq 0 \mathrm{\ when\ }E(0)=0.
\end{eqnarray*}
then the corresponding local solution of the system (\ref{CKGE_P1})
blows up in a finite time $T_{\max}<\infty$, that is,
\begin{eqnarray*}
\lim_{t\rightarrow
T_{\max}^-}(\alpha\|u(t)\|_2^2+\|v(t)\|_2^2)=\infty.
\end{eqnarray*}
\end{theorem}
\textbf{Proof.} We first consider the case $E(0)<0$. The auxiliary
function $G(t)$ in (\ref{G_def}) will still be used here. Naturally
by (\ref{E_conservation}), (\ref{E_def}) and (\ref{G2diff}) we see
that
\begin{eqnarray*}
G^{\prime\prime}(t)&=&(p+q+4)\int
(\alpha|u_t|^2+|v_t|^2)dx+(p+q)\left[\alpha\int(|\nabla
u|^2+m_1^2|u|^2+K_1(x)|u|^2)dx\right. \nonumber\\
&&+\left.\int(|\nabla
v|^2+m_1^2|v|^2+K_2(x)|v|^2)dx\right]-2(p+q+2)E(t)+2b.
\end{eqnarray*}

 Since $E(0)<0$, we now let the constant
$b$ satisfy
\begin{eqnarray*}
0<b\leq -2E(0).
\end{eqnarray*}
Then it follows that
\begin{eqnarray*}
-2(p+q+2)E(t)+2b\geq (p+q+4)b,
\end{eqnarray*}
which implies that
\begin{eqnarray*}
G^{\prime\prime}(t)\geq (p+q+4)\left[\int
(\alpha|u_t|^2+|v_t|^2)dx+b\right].
\end{eqnarray*}

Obviously
\begin{eqnarray}
G^{\prime\prime}(t)>0\label{Th41-G2diff>0}
\end{eqnarray}
on $[0, T_{\max})$.

Moreover we can take a sufficiently large $T_1>0$ and a suitable $b$
such that
\begin{eqnarray}
G^\prime(0)=2\left[\int (\alpha
u_tu+v_tv)dx+bT_1\right]>0\label{Th41-G1diff0>0}
\end{eqnarray}

Thus by (\ref{Th41-G2diff>0}) and (\ref{Th41-G1diff0>0}) we obtain
that $G(t)>0$ and $G^\prime(t)>0$ for every $t\in[0, T_{\max})$.
That is, $G(t)$ and $G^\prime(t)$ are strictly increasing on
$[0,T_{\max})$. Then as in (\ref{G2ttG-G1t>0}) we see that
\begin{eqnarray*}
G^{\prime\prime}(t)G(t)-\frac{p+q+4}{4}(G^\prime(t))\geq 0
\end{eqnarray*}

Thus we have
 \begin{eqnarray}
 \frac{d}{dt}G^{-(p+q)/4}(t)=-\frac{p+q}{4}G^{-(p+q+4)/4}(t)G^\prime(t)<0,\label{G_con1}
 \end{eqnarray}
\begin{eqnarray}
\frac{d^2}{dt^2}G^{-(p+q)/4}(t)&=& -\frac{p+q}{4}G^{-(p+q+8)/4}
\left(G^{\prime\prime}(t)G(t)-\frac{p+q+4}{4}(G^\prime(t))\right)\nonumber\\
&\leq& 0\label{G_con}
\end{eqnarray}
for every $t\in[0, T_{\max})$, which implies that $G^{-(p+q)/4}(t)$
is concave on $[0, T_{\max})$. From (\ref{G_con1}) and (\ref{G_con}) it follows
that the function $G^{-(p+q)/4}\rightarrow 0$ when $t<T_{max}$ and
$t\rightarrow T_{\max}$ ($\displaystyle
T_{\max}\leq\frac{4G(0)}{(p+q)G^\prime(0)}$). Thus we see that there exists
a finite time $T_{\max}>0$ such that
\begin{eqnarray*}
\lim_{t\rightarrow
T_{\max}^-}(\alpha\|u(t)\|_2^2+\|v(t)\|_2^2)=\infty.
\end{eqnarray*}

We next deal with the case $E(0)=0$ with $\displaystyle\int (\alpha
u_0u_1+v_0v_1)dx\geq 0$. Here we define
\begin{eqnarray*}
G(t)=\int(\alpha|u(t)|^2+|v(t)|^2)dx.
\end{eqnarray*}
By direct calculation we have
\begin{eqnarray}
G^\prime(t)=2\int(\alpha u_tu+v_tv)dx
\end{eqnarray}
and
\begin{eqnarray}
G^{\prime\prime}(t)&=&2\int
(\alpha|u_t|^2+|v_t|^2)dx-I(u(t),v(t)).\label{TH41-1}
\end{eqnarray}
By (\ref{E_conservation}) and (\ref{E_def}) we see that
\begin{eqnarray*}
&&\int (\alpha(|\nabla u|^2+m_1^2|u|^2+K_1(x)|u|^2)+(|\nabla
v|^2+m_2^2|v|^2+K_2(x)|v|^2))dx\\
&&\hskip200pt-2a_2^\prime\int|u|^{p+1}|v|^{q+1}dx\leq
0
\end{eqnarray*}
for every $t\in[0, T_{\max})$. Thus we easily have
\begin{eqnarray*}
I(u(t,x),v(x))&=&\int (\alpha(|\nabla
u|^2+m_1^2|u|^2+K_1(x)|u|^2)+(|\nabla
v|^2+m_2^2|v|^2+K_2(x)|v|^2))dx\nonumber\\
&&-(p+q+2)a_2^\prime\int|u|^{p+1}|v|^{q+1}dx<0
\end{eqnarray*}
for every $t\in[0, T_{\max})$.

By (\ref{TH41-1}) we then see that
\begin{eqnarray}
G^{\prime\prime}(t)>0\label{th41-3}
\end{eqnarray}
on $[0, T_{\max})$. And noting the fact $\displaystyle\int(\alpha
u_tu+v_tv)dx\geq 0$, we have
\begin{eqnarray}
G^{\prime}(t)>0\label{TH41-2}
\end{eqnarray}
for every $t\in(0, T_{\max})$.

Thus, by (\ref{th41-3}) and (\ref{TH41-2}) we see that $G(t)$ and
$G^\prime(t)$ are strictly increasing on $[0, T_{\max})$.

Moreover,
\begin{eqnarray*}
G^{\prime\prime}(t)&=&(p+q+4)\int
(\alpha|u_t|^2+|v_t|^2)dx+(p+q)\left[\alpha\int(|\nabla
u|^2+m_1^2|u|^2+K_1(x)|u|^2)dx\right. \nonumber\\
&&+\left.\int(|\nabla
v|^2+m_1^2|v|^2+K_2(x)|v|^2)dx\right]-2(p+q+2)E(0)
\end{eqnarray*}
Noting here $E(0)=0$, we then have
\begin{eqnarray*}
G^{\prime\prime}(t)\geq (p+q+4)\int (\alpha|u_t|^2+|v_t|^2)dx.
\end{eqnarray*}
Since $G(t)>0$ for every $t\in[0, T_{\max})$, we obtain by
Cauchy-Schwartz inequality
\begin{eqnarray*}
G^{\prime\prime}(t)G(t)-\frac{p+q+4}{4}(G^\prime(t))^2\geq 0
\end{eqnarray*}
for every $t\in[0,T_{\max})$.

Then by a concavity argument as in Theorem 4.1, we can claim that there
exists a finite time $T_{\max}<\infty$ such that
\begin{eqnarray*}
\lim_{t\rightarrow T^-_{\max}}(\alpha\|u(t)\|^2+\|v(t)\|^2)=\infty.
\end{eqnarray*}
\begin{flushright}
$\Box$
\end{flushright}
\begin{remark}
For $1<p,q<\infty$ when $n=2$ and for $1<p,q<\frac{2}{n-2}$ when
$n=3$, Theorem 5.1 reproduces the blowing up result in Section 4. But
because of the restriction of the embedding theorem (Lemma
3.1), we cannot apply the method of Section 4 to get the blow up result in this section.
\end{remark}

\section{Blow up with arbitrarily positive initial energy}
To the best of our acknowledge, there is no result for a system of
Klein-Gordon equations when the initial energy is given arbitrarily
positive. In the section we will prove a blow up result for
the system (\ref{CKGE_P}) with arbitrarily positive initial energy.
Indeed, we give the sufficient conditions for the initial datum with
positive initial energy such that the corresponding solution blows
up in a finite time.
\begin{theorem} Assume that $p$ and $q$ satisfy the conditions
(\ref{pq-cond1}), (\ref{pq-cond2}) and (\ref{pq-cond3}). If the
initial data $(u_0, v_0)\in H^1_{0, K_1}\times H^1_{0, K_2}$ and
$(u_1,v_1)\in L^2\times L^2$ satisfy
\begin{eqnarray}
&&E(0)>0;\label{TH51-1}\\
&&I(u_0,v_0)<0;\label{TH51-2}\\
&&\displaystyle\int(\alpha u_0u_1+v_0v_1)dx\geq 0;\label{TH51-3}\\
&&\displaystyle \alpha\|u_0\|^2+\|v_0\|^2>
\frac{2(p+q+2)}{\min\{m^2_1,m^2_2\}(p+q)}E(0).\label{TH51-4}
\end{eqnarray}
Then the corresponding solution $(u(t,x),v(t,x))$ of the system
(\ref{CKGE_P1}) blows up in a finite time $T_{\max}<\infty$, that is,
\begin{eqnarray*}
\lim_{t\rightarrow T_{\max}^-}(\alpha\|u(t)\|^2+\|v(t)\|^2)=\infty.
\end{eqnarray*}
\end{theorem}
\textbf{Proof.} We will prove the result in two steps.

Firstly we show that
\begin{eqnarray}
I(u(t),v(t))<0\label{TH51-5}
\end{eqnarray}
and
\begin{eqnarray}
\alpha\|u(t)\|^2+\|v(t)\|^2>
\frac{2(p+q+2)}{\min\{m^2_1,m^2_2\}(p+q)}E(0)\label{TH51-6}
\end{eqnarray}
for every $t\in[0, T_{\max})$.

We prove (\ref{TH51-5}) by a contradiction argument. Assume that (\ref{TH51-5}) is wrong at some $t\in(0,T_{\max})$, that
is to say, there exists $T>0$ such that
\begin{eqnarray}
T=\min\{t\in[0, T_{\max}); I(u(t),v(t))\geq 0\}.\label{T_def}
\end{eqnarray}
Then by the continuity of $I(u(t),v(t))$ in $t$ we see that
\begin{eqnarray}
I(u(T),v(T))=0.\label{TH51-8}
\end{eqnarray}

Now letting
\begin{eqnarray*}
G(t)=\int(\alpha |u(t,x)|^2+|v(t,x)|^2)dx,
\end{eqnarray*}
we have
\begin{eqnarray}
G^\prime(t)=2\int(\alpha u_tu+v_tv)dx
\end{eqnarray}
and
\begin{eqnarray}
G^{\prime\prime}(t)&=&2\int
(\alpha|u_t|^2+|v_t|^2)dx-I(u(t),v(t)).\label{TH51-10}
\end{eqnarray}
Noting the definition (\ref{T_def}) we see that
\begin{eqnarray}
I(u(t),v(t))<0
\end{eqnarray}
for every $t\in[0, T)$. Thus it follows that $G^{\prime\prime}(t)>0$
on $[0, T)$. And by (\ref{TH51-3}) we have $G^\prime(t)>0$ for
$t\in(0, T)$. In other words, $G(t)$ and $G^\prime(t)$ are
strictly increasing on $[0, T)$. So by (\ref{TH51-4})
\begin{eqnarray}
G(t)>\frac{2(p+q+2)}{\min\{m^2_1,m^2_2\}(p+q)}E(0).\label{TH51-7}
\end{eqnarray}
for every $t\in[0, T)$.

Furthermore, since $u(t)$ and $v(t)$ are continuous in $t$ we get by
(\ref{TH51-7})
\begin{eqnarray}
G(T)>\frac{2(p+q+2)}{\min\{m^2_1,m^2_2\}(p+q)}E(0).\label{Cont-th51-1}
\end{eqnarray}

On the other hand, by (\ref{E_conservation}) and (\ref{E_def}) we
have
\begin{eqnarray*}
\int (\alpha(|\nabla u|^2+m^2_1|u|^2+K_1|u|^2)+(|\nabla
v|^2+m^2_2|v|^2+K_2|v|^2))dx\\
-2a_2^\prime\int|u|^{p+1}|v|^{q+1}dx\leq 2E(0).
\end{eqnarray*}
By (\ref{TH51-8}) we then have
\begin{eqnarray}
G(T)&=&\int (\alpha|u(T)|^2+|v(T)|^2)dx\nonumber\\
&\leq&
\frac{2(p+q+2)}{\min\{m^2_1,m^2_2\}(p+q)}E(0).\label{Cont-th51-2}
\end{eqnarray}

Obviously there is a contradiction between (\ref{Cont-th51-1}) and
(\ref{Cont-th51-2}). Thus we have proved that
\begin{eqnarray}
I(u(t),v(t))<0\label{TH51-9}
\end{eqnarray}
for every $t\in[0, T_{\max})$.

By the argument above we see that $G(t)$ is strictly increasing on
$[0, T_{\max})$ if $I(u(t),v(t))<0$ for every $t\in[0, T_{\max})$
and (\ref{TH51-3}) holds. Namely (\ref{TH51-9}) implies that
\begin{eqnarray}
G(t)>\frac{2(p+q+2)}{\min\{m^2_1,m^2_2\}(p+q)}E(0)\label{TH51-11}
\end{eqnarray}
for every $t\in[0, T_{\max})$.

Now we are going to show the blow up result. By a simple
computation we have
\begin{eqnarray*}
G^{\prime\prime}(t)&=&(p+q+4)\int
(\alpha|u_t|^2+|v_t|^2)dx+(p+q)\left[\alpha\int(|\nabla
u|^2+m_1^2|u|^2+K_1(x)|u|^2)dx\right. \nonumber\\
&&+\left.\int(|\nabla
v|^2+m_2^2|v|^2+K_2(x)|v|^2)dx\right]-2(p+q+2)E(0)\nonumber\\
&\geq& (p+q+4)\int
(\alpha|u_t|^2+|v_t|^2)dx+(p+q)\min\{m^2_1,m^2_2\}\int(\alpha|u|^2+|v|^2)dx\nonumber\\
&&-2(p+q+2)E(0) \nonumber\\
&\geq& (p+q+4)\int (\alpha|u_t|^2+|v_t|^2)dx
\end{eqnarray*}
for every $t\in[0, T_{\max})$.

Thus, by Cauchy-Schwartz inequality we get
\begin{eqnarray*}
G^{\prime\prime}(t)G(t)-\frac{p+q+4}{4}(G^\prime(t))^2&=&(p+q+4)\int
(\alpha|u_t|^2+|v_t|^2)dx\int(\alpha|u|^2+|v|^2)dx\nonumber\\
&&-\int(\alpha u_tu+v_tv)dx\nonumber\\
&\geq& 0.
\end{eqnarray*}

So we have
 \begin{eqnarray}
 \frac{d}{dt}G^{-(p+q)/4}(t)=-\frac{p+q}{4}G^{-(p+q+4)/4}(t)G^\prime(t)<0\label{TH51_con1}
 \end{eqnarray}
\begin{eqnarray}
\frac{d^2}{dt^2}G^{-(p+q)/4}(t)&=& -\frac{p+q}{4}G^{-(p+q+8)/4}
\left(G^{\prime\prime}(t)G(t)-\frac{p+q+4}{4}(G^\prime(t))\right)\nonumber\\
&\leq& 0\label{TH51_con}
\end{eqnarray}
for every $t\in[0, T_{\max})$, which implies that $G^{-(p+q)/4}(t)$
is concave on $[0, T_{\max})$. From (\ref{TH51_con1}) and (\ref{TH51_con}) it follows
that the function $G^{-(p+q)/4}\rightarrow 0$ when $t<T_{max}$ and
$t\rightarrow T_{\max}$ ($\displaystyle
T_{\max}\leq\frac{4G(0)}{(p+q)G^\prime(0)}$). Thus we see that there exists
a finite time $T_{\max}>0$ such that
\begin{eqnarray}
\lim_{t\rightarrow T^-}(\alpha\|u(t)\|_2^2+\|v(t)\|_2^2)=\infty.
\end{eqnarray}
\begin{flushright}
$\Box$
\end{flushright}

\textbf{Acknowledgments.} The author wishes to express his deep
gratitude to Prof. Hitoshi Kitada for his constant encouragement and
careful reading the manuscript. The study is supported by Japanese
Government Scholarship.

\bibliographystyle{amsplain}

\end{document}